\documentclass[12pt]{article}
\usepackage{latexsym,amssymb,amsmath}
\usepackage{amsthm, amstext}
\usepackage{array, amsfonts, mathrsfs}
\usepackage{mathrsfs}
\usepackage[dvips]{graphicx,color}

\newtheorem{Theorem}{Theorem}
\newtheorem{Lemma}{Lemma}

\newtheorem{Proposition}{Proposition}

\def\blue{\color{blue}}
\def\red{\color{red}}

\date{}
\begin{document}
\author{M.I.Belishev\thanks {Saint-Petersburg Department of the Steklov Mathematical Institute, RAS,
                 belishev@pdmi.ras.ru;  Saint-Petersburg State University, 7/9 Universitetskaya nab.,
                 St. Petersburg, 199034, Russia, m.belishev@spbu.ru. Supported by the RFBR grant
                 17-01-00529-à and Volks-Wagen Foundation.} and
        A.F.Vakulenko\thanks{Saint-Petersburg Department of the Steklov Mathematical Institute,
                 vak@pdmi.ras.ru}}
\title{On algebras of harmonic quaternion fields in ${\mathbb R}^3$}
\maketitle

\begin{abstract}
Let ${\mathscr A}(D)$ be an algebra of functions continuous in the
disk $D=\{z\in{\mathbb C}\,|\,\,\,|z|\leqslant 1\}$ and {\it
holomorphic} into $D$. The well-known fact is that the set
${\mathscr M}$ of its characters (homomorphisms ${\mathscr
A}(D)\to\mathbb C$) is exhausted by the Dirac measures
$\{\delta_{z_0}\,|\,\,z_0\in D\}$ and a homeomorphism ${\mathscr
M}\cong D$ holds. We present a 3d analog of this classical result
as follows.

Let $B=\{x\in{\mathbb R}^3\,|\,\,|x|\leqslant 1\}$. A quaternion
field is a pair $p=\{\alpha,u\}$ of a function $\alpha$ and vector
field $u$ in the ball $B$. A field $p$ is {\it harmonic} if
$\alpha, u$ are continuous in $B$ and $\nabla\alpha={\rm
rot\,}u,\,{\rm div\,}u=0$ holds into $B$. The space ${\mathscr
Q}(B)$ of such fields is not an algebra w.r.t. the relevant
(point-wise quaternion) multiplication. However, it contains the
commutative algebras ${\mathscr A}_\omega(B)=\{p\in{\mathscr
Q}(B)\,|\,\,\nabla_\omega\alpha=0,\,\nabla_\omega
u=0\}\,\,(\omega\in S^2)$, each ${\mathscr A}_\omega(B)$ being
isometrically isomorphic to ${\mathscr A}(D)$. This enables one to
introduce a set ${\mathscr M}^{\mathbb H}$ of the $\mathbb
H$-valued linear functionals on ${\mathscr Q}(B)$ ({\it $\mathbb
H$-characters}), which are multiplicative on each ${\mathscr
A}_\omega(B)$, and prove that ${\mathscr M}^{\mathbb
H}=\{\delta^{\mathbb H}_{x_0}\,|\,\,x_0\in B\}\cong B$, where
$\delta^{\mathbb H}_{x_0}(p)=p(x_0)$.

\end{abstract}


\bigskip

\setcounter{section}{-1}

\section{Introduction}\label{sec Introduction}

\noindent{$\bullet$}\,\,\,The result, which our paper is devoted
to, is announced in \cite{B Quat 2016}\footnote{unfortunately,
with some inaccuracies in the formulations}. In a sense, it is a
`by-product' of activity in the framework of the algebraic
approach to tomography problems on manifolds \cite{B Calderon
2003}--\cite{BSharaf 2008}. However, we hope that this result is
of certain independent interest for the real uniform Banach
algebras theory \cite{Abel Jarosz, Jarosz}. Namely, we propose a
3d generalization of the well-known theorem on the characters of
the disk-algebra of holomorphic functions. Perhaps, the most
curious point is that such a generalization does exist although
the relevant 3d analog of the disc-algebra {\it is not an
algebra}.
\smallskip

\noindent{$\bullet$}\,\,\,A subject to be generalized is the
following well-known result.

Let ${\mathscr A}(D)$ be a commutative Banach algebra of functions
continuous in the disk $D=\{z\in{\mathbb C}\,|\,\,\,|z|\leqslant
1\}$ and holomorphic into $D$. The well-known fact (see, e.g.,
\cite{Naimark}) is that the set ${\mathscr M}$ of its {\it
characters} (homeomorphisms ${\mathscr A}(D)\to\mathbb C$) endowed
with the Gelfand topology, is exhausted by the Dirac measures
$\{\delta_{z_0}\,|\,\,z_0\in D\}$ and the homeomorphism ${\mathscr
M}\cong D$ holds.
\smallskip

Shortly, a 3d generalization, which is our main result, looks as
follows.

Let $B=\{x\in{\mathbb R}^3\,|\,\,|x|\leqslant 1\}$. A quaternion
field is a pair $p=\{\alpha,u\}$ of a function $\alpha$ and vector
field $u$ in the ball $B$. Such fields are identified with the
${\mathbb H}$-valued functions and, hence, can be multiplied
point-wise as quaternions. We say a field $p$ to be {\it harmonic}
if $\alpha, u$ are continuous in $B$ and $\nabla\alpha={\rm
rot\,}u,\,{\rm div\,}u=0$ holds into $B$. A space ${\mathscr
Q}(B)$ of such fields is not an algebra w.r.t. the above mentioned
multiplication. However, it contains the commutative algebras
${\mathscr A}_\omega(B)=\{p\in{\mathscr
Q}(B)\,|\,\,\nabla_\omega\alpha=0,\,\nabla_\omega
u=0\}\,\,(\omega\in S^2)$, which we call the {\it axial algebras}
($\omega$ is an axis). Each ${\mathscr A}_\omega(B)$ is
isometrically isomorphic to ${\mathscr A}(D)$. This enables one to
introduce the set ${\mathscr M}^{\mathbb H}$ of the $\mathbb
H$-valued linear functionals on ${\mathscr Q}(B)$ ({\it $\mathbb
H$-characters}), which are multiplicative on all ${\mathscr
A}_\omega(B)$, and prove that ${\mathscr M}^{\mathbb
H}=\{\delta^{\mathbb H}_{x_0}\,|\,\,x_0\in B\}\cong B$, where
$\delta^{\mathbb H}_{x_0}(p)=p(x_0)$.

Sections  \ref{sec The classics},\ref{sec The 3d analog} contain
the definitions, statements of the results and comments on them.
The basic proofs are placed in \ref{sec Proofs}.
\smallskip

\noindent{$\bullet$}\,\,\,We'd like to thank Dr C.Shonkwiler for
helpful remarks and very useful references.

\section{The 2d classics}\label{sec The classics}
We begin with the result, which is planned to be generalized.

\subsubsection*{Algebras}
\noindent$\bullet$\,\,\,Let $D=\{z\in{\mathbb
C}\,|\,\,\,|z|\leqslant 1\}$ be the disc on the complex plane. The
continuous function algebra
 \begin{align*}
& C^{\mathbb C}(D)\,:=\,\{f=\varphi+i\psi\,|\,\,\,\varphi,\psi \in
C^{\mathbb
R}(D)\};\quad \|f\|\,=\,\underset{D}{\rm sup}\,|f|\,,\\
& \|fg\|\leqslant\|f\|\|g\|\,,\quad fg=gf\,,\quad \|f^2\|=\|f\|^2
 \end{align*}
is a Banach commutative uniform algebra.
\smallskip

\noindent$\bullet$\,\,\, A function $f=\varphi+i\psi$ is {\it
holomorphic} if the Cauchy-Riemann conditions
 \begin{equation}\label{Eq CR 1}
d\psi\,=\,\star\,d\varphi \qquad (\delta\psi=\delta\varphi=0)
 \end{equation}
hold in the inner points of ${\rm Dom\,}f\subset\mathbb C$. Here
$\varphi$ and $\psi$ are regarded as the 0-forms, $\star$ is the
Hodge operator corresponding to the standard orientation of
$\mathbb C$, $d$ and $\delta$ are the differential and
codifferential respectively. The conditions in the brackets are
fulfilled just by the well-known definitions. It is a form of
writing the CR-conditions, which is most relevant to the
forthcoming generalization.
\smallskip

\noindent$\bullet$\,\,\,The (sub)algebra
 \begin{align*}
{\mathscr A}(D)\,:=\,\{f\in C^{\mathbb C}(D)\,|\,\,\,f\,\,\text{is
holomorphic into}\,\,D\}
 \end{align*}
is also a commutative uniform Banach algebra.

\subsubsection*{Characters}
By ${\mathfrak L}(F,G)$ we denote the normed space of the linear
continuous operators from a Banach space $F$ to a Banach space
$G$.

Let ${\mathscr A}^\prime(D):={\mathfrak L}({\mathscr
A}(D),{\mathbb C})$ be the dual space. {\it Characters}
(multiplicative functionals) are defined as elements of the set
 \begin{align*}
{\mathscr M}\,:=\,\left\{\mu \in {\mathscr
A}^\prime(D)\,|\,\,\,\mu(fg)=\mu(f)\mu(g),\,\,\,\,\,f,g\in
{\mathscr A}(D)\right\}
 \end{align*}
endowed with the Gelfand ($\ast$-weak) topology, which is
determined by the convergence
 \begin{align*}
\{\mu_j\to \mu\}\Leftrightarrow\{\mu_j(f)\overset{\mathbb C}\to
\mu(f)\,,\,\,\,f\in {\mathscr A}(D)\}\,.
 \end{align*}
The set ${\mathscr M}$ is also called a {\it spectrum} of the
algebra ${\mathscr A}(D)$.

An example of characters, which turns out to be universal, is
provided by the {\it Dirac measures} \,$\delta_{z_0}\in {\mathscr
M}$:
 $$
\delta_{z_0}(f)\,:=\,f(z_0)\,, \qquad f\in{\mathscr
A}(D)\quad\,\,(z_0\in D)\,.
 $$

\subsubsection*{Basic fact}
For topological spaces, we wright $S\cong T$ if $S$ and $T$ are
homeomorphic. For algebras, ${\mathscr A}\cong {\mathscr B}$ means
that ${\mathscr A}$ and ${\mathscr B}$ are isometrically
isomorphic.

Our goal is to provide a 3d analog of the following classical
result (see, e.g., \cite{Naimark}: Chapter III, paragraph 11, item
3).
\begin{Theorem}\label{T1}
For any $\mu\in{\mathscr M}$ there is a point $z^\mu\in D$ such
that $\mu=\delta_{z^\mu}$. The map
 $$
\gamma: {\mathscr M}\to D,\quad \mu\,\mapsto\,z^\mu
 $$
is a homeomorphism, so that ${\mathscr M}\cong D$ holds. The
Gelfand transform
 $$
\Gamma: {\mathscr A}(D)\to C^{\mathbb C}({\mathscr
M}),\quad(\Gamma f)(\mu)=\mu(f),\,\,\,\mu\in {\mathscr M}
 $$
is an isometric isomorphism onto its image, so that ${\mathscr
A}(D)\cong \Gamma\left({\mathscr A}(D)\right)$ holds.
\end{Theorem}
Thus, the spectrum of ${\mathscr A}(D)$ is exhausted by Dirac
measures: ${\mathscr M}=\{\delta_{z_0}\,|\,\,z_0\in D\}$. Also,
since $(\Gamma f)(\mu)=\mu(f)=\delta_{z^\mu}(f)=f(z^\mu)$, the
Gelfand transform just transfers functions from $D$ to ${\mathscr
M}$ along the map $\gamma$.
\smallskip

Notice in addition that the same set of the Dirac measures
exhausts the spectrum of the `big' algebra $C^{\mathbb C}(D)$
\cite{Naimark}.

\section{The 3d analogs}\label{sec The 3d analog}

\subsubsection*{Quaternions}
\noindent$\bullet$\,\,\,Recall that ${\mathbb H}$ is a real
algebra of the collections ({\it quaternions})
 $$
{\mathfrak h}=\alpha +u_1{\bf i}+u_2{\bf j}+u_3{\bf
k},\,\qquad\alpha, u_i\in{\mathbb R}
 $$
endowed with the component-wise linear operations and the norm
(module) $|{\mathfrak
h}|=[\alpha^2+u_1^2+u_2^2+u_3^2]^{\frac{1}{2}}$. A multiplication
is determined by the table
$$
{\bf i}{\bf i}={\bf j}{\bf j}={\bf k}{\bf k}=-1;\quad {\bf i}{\bf
j}={\bf k},\,\,{\bf j}{\bf k}={\bf i},\,\,{\bf k}{\bf i}={\bf j}
$$
and extended to ${\mathbb H}$ by linearity and distributivity. The
multiplication is associative but not commutative. The module
obeys $|{\mathfrak g}{\mathfrak h}|=|{\mathfrak g}|\,|{\mathfrak
h}|$.
\smallskip

\noindent$\bullet$\,\,\,A {\it geometric quaternion} is a pair of
a scalar and 3d vector. The set of pairs
 $$
{\mathbb H}_{\rm g}=\left\{p=\{\alpha,
u\}\,|\,\,\,\alpha\in{\mathbb R},\,\,u\in{\mathbb R}^3\right\}
 $$
is a real algebra w.r.t. the component-wise linear operations, the
module $|p|=[\alpha^2+|u|^2]^{\frac{1}{2}}$, and multiplication
 $$
pq\,:=\,\{\alpha\beta-u\cdot v,\,\,\alpha v+\beta u+u\wedge
v\}\quad \text{for}\,\,p=\{\alpha, u\},\,\,q=\{\beta,v\}\,,
 $$
where $\cdot$ and $\wedge$ are the standard inner and vector
products in ${\mathbb R}^3$. The multiplication is noncommutative.
The module obeys $|pq|=|p||q|$.

The correspondence
 \begin{equation}\label{Eq quaternion-vector}
{\mathbb H}\ni {\mathfrak h}=\alpha +u_1{\bf i}+u_2{\bf j}+u_3{\bf
k}\,\leftrightarrow\, h=\{\alpha, \begin{pmatrix} u_1\cr u_2\cr
u_3\end{pmatrix}\}\in{\mathbb H}_{\rm g}
 \end{equation}
determines an isometric isomorphism of algebras: ${\mathbb
H}\cong{\mathbb H}_{\rm g}$. As a consequence, ${\mathbb H}_{\rm
g}$ is the left ${\mathbb H}$-module with the action
 \begin{equation}\label{Eq H-module}
{\mathfrak h}p\,=\,hp\,, \qquad {\mathfrak h}\in {\mathbb
H},\,p\in{\mathbb H}_{\rm g}\,.
 \end{equation}

\noindent$\bullet$\,\,\,A {\it quaternion field} is an ${\mathbb
H}_{\rm g}$-valued function on a domain in ${\mathbb R}^3$, i.e.,
a pair $p=\{\alpha, u\}$\,, where $\alpha$ is a function and $u$
is a vector field defined on a domain in ${\mathbb R}^3$. The set
of such fields is a real algebra and left ${\mathbb H}$-module
w.r.t. the relevant point-wise operations.

\subsubsection*{Algebras and spaces}
\noindent$\bullet$\,\,\,Let $B=\{x\in{\mathbb
R}^3\,|\,\,\,|x|\leqslant 1\}$ be a ball. The space of continuous
fields
 \begin{align*}
C^{\mathbb H}(B)\,:=\,\{\,p=\{\alpha,u\}\,|\,\,\alpha \in
C^{\mathbb R}(B),\,\,u \in C(B; {\mathbb R}^3)\}
 \end{align*}
with the norm $\|p\|\,=\,\underset{B}{\rm sup}\,|p|$ obeying
$\|pq\|\leqslant\|p\|\|q\|$ and $\|p^2\|=\|p\|^2$ is a
(noncommutative) Banach uniform algebra and a left ${\mathbb
H}$-module.
\smallskip

\noindent$\bullet$\,\,\,A field $p=\{\alpha,u\}$ is said to be
{\it harmonic} if the Cauchy-Riemann conditions
 \begin{equation}\label{Eq CR 2}
d\alpha\,=\,\star\,du',\quad \delta u'=0 \qquad(\delta\alpha=0)
 \end{equation}
hold in the inner points of ${\rm Dom\,}p\subset{\mathbb R}^3$.
Here $\alpha$ is regarded as a 0-form, $\star$ is the Hodge
operator corresponding to the standard orientation of ${\mathbb
R}^3$, $d$ and $\delta$ are a differential and codifferential
respectively, $u'$ is a 1-form dual to $u$, i.e., $u'(b)=b\cdot u$
on vector-fields $b$. The condition in the brackets is fulfilled
automatically, whereas the `gauge condition' $\delta u'=0$ is now
not trivial. In terms of the vector analysis operations, (\ref{Eq
CR 2}) is equivalent to
 \begin{equation}\label{Eq CR 3}
\nabla \alpha\,=\,{\rm rot\,}u\,,\quad{\rm div\,}u=0 \qquad {\rm
into\,\,Dom\,}p\,.
 \end{equation}

\noindent$\bullet$\,\,\,The key object of the paper is the
(sub)space of harmonic fields
 \begin{align*}
{\mathscr Q}(B)\,=\,\{p\in C^{\mathbb H}(B) \,|\,\,p\,\, \text{is
harmonic into}\,\,B\}\,.
 \end{align*}
A simple calculation with regard to the definitions (\ref{Eq
H-module}) and (\ref{Eq CR 3}) enables one to check the following
property (see \cite{B Quat 2016} for detail).
\begin{Proposition}\label{Prop Q(B) left module}
The space ${\mathscr Q}(B)$ is a left $\mathbb H$-module:
 $
{\mathfrak h} \in {\mathbb H}\,\,\,\text{and}\,\,\,p\in {\mathscr
Q}(B)\,\,\,\text{imply}\quad {\mathfrak h}p\in {\mathscr Q}(B)\,.
 $
 \end{Proposition}
\noindent In the mean time, ${\mathscr Q}(B)$ is not a
(sub)algebra:\, generically, $p,q\in{\mathscr Q}(B)$ doesn't imply
$pq\in{\mathscr Q}(B)$. It is the fact, which was perceived as an
obstacle for existence of a 3d version of Theorem \ref{T1}.
However, we'll see that such a version does exist, whereas the
space ${\mathscr Q}(B)$ turns out to be a relevant analog of the
algebra ${\mathscr A}(D)$.

A constant quaternion field, which is equal to $h\in{\mathbb H}_g$
identically, will be denoted by the same symbol $h$. Such fields
belong to ${\mathscr Q}(B)$.
\smallskip

\noindent$\bullet$\,\,\,Let $\omega\in S^2$ be a unit vector.

For a function $\alpha$, we denote
$\nabla_\omega\alpha=\omega\cdot\nabla\alpha$ and say $\alpha$ to
be {\it $\omega$-axial} if $\nabla_\omega\alpha=0$. For a vector
field $u$, we denote by $\nabla_\omega u$ its covariant derivative
in ${\mathbb R}^3$ and say $u$ to be {\it $\omega$-axial} if
$\nabla_\omega u=0$. A quaternion field $p=\{\alpha,u\}$ is
$\omega$-axial if its components are $\omega$-axial. To be
$\omega$-axial just means to be constant on the straight lines
$x=x_0+t\omega\,\,\,\,(x_0\in{\mathbb R}^3,\,\,t\in \mathbb R)$.
\smallskip

\noindent$\bullet$\,\,\,The following is proven in \cite{B Quat
2016}.
 \begin{Proposition}\label{Prop nabla psi=omega vedge nabla phi}
An axial field $p=\{\varphi, \psi\omega\}$ is harmonic if and only
if the functions $\varphi,\,\psi$ obey
 \begin{equation*}
\nabla\psi\,=\,\omega\wedge\nabla\varphi\qquad {\rm
into\,\,\,Dom\,}p\,.
 \end{equation*}
In this case, $\varphi$ and $\psi$ are $\omega$-axial and
harmonic: $\Delta\varphi=\Delta\psi=0$ holds into ${\rm Dom\,}p$.
\end{Proposition}

If $q=\{\lambda, \rho\omega\}$, then
 $$
pq\,=\,\{\varphi\lambda-\psi\rho,\,[\varphi\rho+\psi\lambda]\,\omega\}\,,
 $$
i.e., $pq$ is also $\omega$-axial and $pq=qp$. Thus, one can
multiply coaxial fields, the multiplication being commutative.
\smallskip

\noindent$\bullet$\,\,\,By the aforesaid, the subspace
 $$
{\mathscr A}_\omega(B)\,=\,\{p\in {\mathscr Q}(B)\,|\,\,p\,\,\,
\text{is}\,\,\omega\text{-axial}\,\}\subset{\mathscr Q}(B)
\qquad\quad(\omega\in S^2)
 $$
is a Banach {\it commutative} uniform algebra. Moreover, each
${\mathscr A}_\omega(B)$ is isometric to ${\mathscr A}(D)$ via the
map $p\mapsto f_p=\tilde\varphi+i\tilde\psi$, where
$\tilde\varphi=\varphi|_{D_\omega}, \tilde\psi=\psi|_{D_\omega}$,
and $D_\omega=\{x\in B\,|\,\,x\cdot\omega=0\}$ is the disc
properly oriented and identified with $D\subset\mathbb C$.

Thus, being not an algebra, the harmonic space contains algebras.
A reserve of these algebras is rich enough: the following fact
will be established later in sec \ref{sec Proofs}.
 \begin{Lemma}\label{Lemma on density}
The relation
  \begin{equation}\label{Eq span A omega = Q(B)}
\overline{{\rm span\,}\left\{{\mathscr
A}_\omega(B)\,|\,\,\omega\in S^2\right\}}\,=\,{\mathscr Q}(B)
  \end{equation}
is valid\,\,\,(the closure in $C^{\mathbb H}(B)$).
 \end{Lemma}

\subsubsection*{$\mathbb H$-characters}
Recall that ${\mathfrak L}(F,G)$ is the space of linear operators
from $F$ to $G$.
\smallskip

\noindent$\bullet$\,\,\,In the 3d case, a role of the dual space
${\mathscr A}^\prime(D)$ is played by the space
 $$
{\mathscr Q}^\times(B)\,=\,\{l\in {\mathfrak L}\left({\mathscr
Q}(B),{\mathbb H}\right)|\,\,\,l({\mathfrak h}p)={\mathfrak
h}l(p),\,\,\,\, \forall p\in {\mathscr Q}(B),{\mathfrak
h}\in\mathbb H\}\,,
 $$
which we call an  {\it $\mathbb H$-dual} to ${\mathscr Q}(B)$; its
elements are named by {\it $\mathbb H$-functionals}. By this
definition, one has ${\mathfrak h}l(p)=l({\mathfrak
h}p)\overset{(\ref{Eq quaternion-vector})}=l(hp)=l(h)l(p)$, which
implies
 \begin{equation}\label{Eq l(h)=h}
l(h)\,=\,{\mathfrak h}\,,\qquad {\mathfrak h}\in \mathbb H
  \end{equation}
for all constant fields $h$.

In addition, note that ${\mathscr Q}^\times(B)$ can be endowed
with a left $\mathbb H$-module structure \cite{Joyce}.
\smallskip

\noindent$\bullet$\,\,\,In the capacity of a 3d-analog of the
Gelfand spectrum $\mathscr M$ of algebra ${\mathscr A}(D)$, we
propose the set
 \begin{align*}
{\mathscr M}^{\mathbb H}\,:=\,\left\{\mu \in {\mathscr
Q}^\times(B)\,|\,\,\,\mu(pq)=\mu(p)\mu(q),\,\,\,\,\,\forall\,p,q\in
{\mathscr A}_\omega(B),\,\omega\in S^2\right\}
 \end{align*}
endowed with $*$-weak topology determined by the convergence
 \begin{align*}
\{\mu_j\to \mu\}\Leftrightarrow\{\mu_j(f)\overset{\mathbb H}\to
\mu(f)\,,\,\,\, \forall f\in {\mathscr Q}(B)\}\,.
 \end{align*}
It looks reasonable to call ${\mathscr M}^{\mathbb H}$ an {\it
$\mathbb H$-spectrum} of the harmonic space ${\mathscr Q}(B)$ and
name its elements by {\it $\mathbb H$-characters}.

An example of $\mathbb H$-characters is provided by the
`quaternion Dirac measures'
 $$
\delta^{\mathbb H}_{x_0}(p)\,=\,p(x_0)\,, \qquad p\in{\mathscr
Q}(B)\quad\,\,(x_0\in B)\,.
 $$
Note that $\delta^{\mathbb H}_{x_0}$ is well defined and {\it
multiplicative} on the `big' algebra $C^{\mathbb H}(B)$.

\subsubsection*{Main result}
The above mentioned example turns out to be universal: as will be
proven later, the $\mathbb H$-spectrum is exhausted by the
quaternion Dirac measures.
 \begin{Theorem}\label{T2}
For any $\mu\in{\mathscr M}^{\mathbb H}$, there is a point
$x^\mu\in B$ such that $\mu=\delta^{\mathbb H}_{x^\mu}$. The map
 $$
\gamma: {\mathscr M}^{\mathbb H}\to B,\quad \mu\,\mapsto\,x^\mu
 $$
is a homeomorphism, so that ${\mathscr M}^{\mathbb H}\cong B$
holds. The Gelfand transform
 $$
\Gamma: {\mathscr Q}(B)\to C({\mathscr M^{\mathbb H}};{\mathbb
H}),\quad(\Gamma f)(\mu):=\mu(f),\,\,\,\mu\in {\mathscr
M}^{\mathbb H}
 $$
is an isometry onto its image, so that ${\mathscr Q}(B)\cong
\Gamma\left({\mathscr Q}(B)\right)$ holds.
 \end{Theorem}
Thus, $\Gamma$ just transfers the fields from $B$ to ${\mathscr
M}^{\mathbb H}$ along the map $\gamma$.
\smallskip

Notice in addition that Dirac measures exhaust the set of $\mathbb
H$-linear functionals of the `big' algebra $C^{\mathbb H}(B)$
\cite{Abel Jarosz, Jarosz}.

\subsubsection*{Summary}
For reader's convenience, we present a correspondence table
between the classical objects and their 3d-analogs.
 \begin{align*}
& {\blue\text{disk}\,D\subset{\mathbb C}}  && {\red\text{ball}\,B\subset{\mathbb R}^3}\\
& {\blue\text{algebra}\,\,C^{\mathbb C}(D)} &&
{\red\text{algebra}\,\,C^{\mathbb H}(B)}\\
& {\blue f=\varphi+i\psi\!:\,\,\,d\psi=\star\,d\varphi} && {\red
p=\{\alpha,u\}\!:\,\,\,du'=\star\,d\alpha,
\,\,\, \delta u'=0}\\
& {\blue \text{holom.\,func.\! algebra}\,\,{\mathscr A}(D)}
&&{\red
\text{harm.\,quat.\,field\! space}\,\,{\mathscr Q}(B)}\\
& {\blue \dots\dots\dots\dots\dots\dots\dots\dots}  && {\red
\text{axial algebras}\,\,{\mathscr A}_\omega(B)\subset{\mathscr Q}(B)}\\
& {\blue {\mathbb C}\text{-linear functionals}\,\, {\mathscr
A}^\prime(D)} && {\red {\mathbb H}\text{-linear functionals}\,\,
{\mathscr Q}^\times(B)}\\
& {\blue \text{characters}\,\, {\mathscr M}\subset{\mathscr
A}^\prime(D)} && {\red {\mathbb H}\text{-characters}\,\, {\mathscr
M}^{\mathbb
H}\subset{\mathscr Q}^\times(B)}\\
& {\blue{\mathscr M}\cong D} && {\red {\mathscr M}^{\mathbb
H}\cong B}
\end{align*}

\subsubsection*{Comments}

\noindent$\bullet$\,\,\,As is mentioned in Introduction, our
results are obtained in the framework of algebraic version of the
so-called BC-method, which is an approach to inverse problems of
mathematical physics \cite{B Calderon 2003}, \cite{B Sobolev Geom
Rings 2008}--\cite{BSharaf 2008}. More precisely, the impact comes
from the impedance tomography of 3d Riemannian manifolds \cite{B
CUBO 3d tomogr 2005, B UMN 2017}. To answer the following
questions would be helpful for the progress in this application.
The questions are put not for a ball $B$ but a domain
$\Omega\subset{\mathbb R}^3$ that should not lead to confusion
since the proper generalizations are evident. However, we don't
know the answers even for a ball.

$1.$\,\,\,For an algebra $\mathscr A$ and a set $S\subset \mathscr
A$, by $\vee S$ we denote a minimal (sub)algebra in $\mathscr A$,
which contains $S$. Does
 $$
\overline{\vee {\mathscr Q}(\Omega)}\,=\,C^{\mathbb H}(\Omega)
 $$
hold at least for a class of $\Omega$'s? Presumably, the answer
may be got by the proper application of Corollary 1 from
\cite{Jarosz}.

$2.$\,\,\,For $p\in C^{\mathbb H}(\Omega)$, we denote
$p^\partial=p|_{\partial\Omega}$. A simple fact is that a harmonic
quaternion field is determined by its boundary values
\footnote{Moreover, owing to the relevant maximal principle, the
map ${\mathscr Q}(\Omega)\ni p\mapsto p^\partial\in C^{\mathbb
H}(\partial\Omega)$ preserves the norms \cite{B Quat 2016}.}.
Therefore, the algebras $\vee{\mathscr
Q}(\Omega)\,\subset\,C^{\mathbb H}(\Omega)$ and
$\vee\{p^\partial\,|\,\,p\in {\mathscr
Q}(\Omega)\}\,\subset\,C^{\mathbb H}(\partial\Omega)$ turn out to
be isomorphic (but not isometric!).

A {\it boundary algebra}
 $$
{\mathscr B}=\overline{\vee\{p^\partial\,|\,\,p\in {\mathscr
Q}(\Omega)\}}
 $$
(the closure in $C^{\mathbb H}(\partial\Omega)$) is a
noncommutative Banach uniform algebra \footnote{The 2d version of
algebra $\mathscr B$ solves the 2d tomography problem on
manifolds: see \cite{B Calderon 2003}, \cite{B Sobolev Geom Rings
2008}, \cite{B UCLA 2013}}. What is the relation between
${\mathscr B}$ and $C^{\mathbb H}(\Omega)$? Let ${\mathfrak
M}({\mathscr B})$ be a structure space of ${\mathscr B}$
\cite{Naimark}; is there a chance for ${\mathfrak M}({\mathscr
B})\cong \Omega$?
\smallskip

In fact, all of these questions are of auxiliary character. They
become full-valued and important if $\Omega$ is a Riemannian
manifold with boundary. In such a case, the algebra ${\mathscr B}$
is also well defined and seems to be a most promising device for
solving 3d impedance tomography problem: to recover $\Omega$ via
its Dirichlet-to-Neumann operator \cite{BSharaf 2008, B UMN 2017}.
Any informative results on ${\mathscr B}$ are welcomed.
\smallskip

\noindent$\bullet$\,\,\,Theorem \ref{T1} is valid not only for a
disc but a much wider class of domains on $\mathbb C$ (see, e.g.,
\cite{Stout}). In the mean time, one can show that Theorem
\ref{T2} remains true for a convex bounded $\Omega\subset{\mathbb
R}^3$. Also, it is valid for a torus $\{x\in{\mathbb
R}^3\,|\,\,{\rm dist\,}(x,L)\leqslant\varkappa R\}$, where
$L\subset{\mathbb R}^3$ is a circle of radius $R$ and
$0<\varkappa<1$. However, the class of available $\Omega$'s is not
properly specified yet.
\smallskip

\noindent$\bullet$\,\,\,Our considerations enable one to set up a
`3d corona problem' for the space of harmonic quaternion fields
{\it bounded} in $B$ (by analogy with $H^\infty(D)$).
\smallskip

It would be interesting to extend Theorem \ref{T2} to
$\Omega\in{\mathbb R}^n$. Presumably, such an extension has to
deal with the spaces and algebras of harmonic differential forms
\cite{BSharaf 2008}.

\section{Proofs}\label{sec Proofs}

In what follows, identifying the quaternions ${\bf i,\,j,\,k}$
with elements of the standard basis in ${\mathbb R}^3$ (see
(\ref{Eq quaternion-vector})), we write $x=x_1{\bf i}+x_2{\bf
j}+x_3{\bf k}\in {\mathbb R}^3$.

Recall that for the vector fields $u=u_1{\bf i}+u_2{\bf j}+u_3{\bf
k}$ one defines the Laplacian by $\Delta u=\Delta u_1\,{\bf
i}+\Delta u_2\,{\bf j}+\Delta u_3\,{\bf k}$. The harmonic vector
fields are the ones satisfying $\Delta u=0$.

\subsubsection*{Polynomials}
Let

$\Pi$ be a linear space of the (scalar) polynomials of the
variables $x_1,x_2,x_3$; $\Pi_n\subset \Pi$ the polynomials of
degree $n\geqslant 0$; $\dot \Pi_n\subset \Pi_n$ the homogeneous
polynomials, i.e., a linear span of monomials
$x_1^{r_1}x_2^{r_2}x_3^{r_3}$ with $r_1+r_2+r_3=n$;

$P=\{\alpha\in \Pi\,|\,\,\Delta \alpha=0\},\, P_n=\{\alpha\in
\Pi_n\,|\,\,\Delta \alpha=0\},\,\dot P_n=\{\alpha\in \dot
\Pi_n\,|\,\,\Delta \alpha=0\}$ the harmonic polynomials;

${\mathscr P}=\{p=\{\alpha,u\}\,|\,\,u=u_1{\bf i}+u_2{\bf
j}+u_3{\bf k}:\,\,\alpha, u_k\in \Pi;\,\,\nabla\alpha={\rm
rot\,}u,\,\,{\rm div\,}u=0\}$ the harmonic quaternion polynomials;
the harmonicity implies $\alpha, u_k \in P$;

${\mathscr P}_n=\{\,\{\alpha,u\}\in {\mathscr P}\,|\,\,u=u_1{\bf
i}+u_2{\bf j}+u_3{\bf k}:\,\,\alpha, u_k\in P_n\}$ the harmonic
quaternion polynomials of degree $n$;

$\dot{\mathscr P}_n=\{\,\{\alpha,u\}\in {\mathscr
P}\,|\,\,u=u_1{\bf i}+u_2{\bf j}+u_3{\bf k}:\,\,\alpha, u_k\in
\dot P_n\}$ the homogeneous harmonic quaternion polynomials of
degree $n$;
\medskip

The notations $\{...\}^\omega$ mean that a set $\{...\}$ consists
of the objects (functions, vector fields, quaternion fields, etc),
which are $\omega$-axial, i.e., take constant values on the
straight lines $\{x=x_0+t\omega\,|\,\,x_0\in{\mathbb
R}^3,\,t\in\mathbb R\}$. So, $P^\omega, \Pi^\omega, {\mathscr
P}^\omega, \dots$ are of clear meaning.
\begin{Lemma}\label{Lemma P=sum Pn homog}
The relation
 \begin{equation}\label{Eq P=sum Pn homog}
{\mathscr P}={\rm span\,}\{{\mathscr P}^\omega\,|\,\,\omega\in
S^2\}
 \end{equation}
holds.
\end{Lemma}
{\bf Proof.}

\noindent$\bullet$\,\,\,The relation $(x_1+x_2{\bf
i})^n=R_n(x_1,x_2)+I_n(x_1,x_2)\bf i$ determines the $\bf k$-axial
polynomials $R_n,I_n\in \dot P_n$ satisfying $\nabla I_n={\bf
k}\wedge\nabla R_n$. Any harmonic homogeneous polynomial of
variables $x_1,\,x_2$ of degree $n$ is $a R_n(x_1,x_2)+b
I_n(x_1,x_2)$\,with $a,b\in \mathbb R$. Hence, ${\rm dim\,}\dot
P_n^{\bf k}=2$. Quite analogously, one has
 \begin{equation}\label{Eq dim dotPn=2}
{\rm dim\,}\dot P_n^{\omega}=2\qquad\text{for any}\,\,\, \omega\in
S^2\,.
\end{equation}
\noindent$\bullet$\,\,\,We omit the proof of the following result,
which is simply derived by induction.
 \begin{Proposition}\label{Prop dim span}
Let $\omega_1,\,\dots ,\,\omega_r \in S^2$ be the pair-wise
different vectors: $\omega_i\not=\pm\omega_j$. The relation
 \begin{equation}\label{Eq dim span dotPn}
{\rm dim\,\,span\,}\left\{\dot P^{\omega_k}_n\,\big|\,\,k=1,
\dots, r\right\}=
                      \begin{cases}
             2r,    & r\leqslant n \\
             2n+1,  & r>n
                      \end{cases}
 \end{equation}
holds.
 \end{Proposition}
\smallskip

\noindent$\bullet$\,\,\,Let $\omega_1,\,\dots ,\,\omega_{n+1} \in
S^2$, $\omega_i\not=\pm\omega_j$. By relations (\ref{Eq dim
dotPn=2}) and (\ref{Eq dim span dotPn}), the subspaces $\dot
P_n^{\omega_1},\dots, \dot P_n^{\omega_{n}}\subset \dot P_n$ are
linearly independent, whereas $\dot P_n^{\omega_1},\dots, \dot
P_n^{\omega_{n+1}}$ are not independent, and ${\rm dim\,}\dot
P_n^{\omega_{n+1}}\cap{\rm span\,}\{\dot P_n^{\omega_1},\dots,
\dot P_n^{\omega_{n}}\}=1$ . Therefore, choosing a nonzero
$\varphi_{n+1}\in\dot P_n^{\omega_{n+1}}\cap{\rm span\,}\{\dot
P_n^{\omega_1},\dots, \dot P_n^{\omega_{n}}\}$, one determines the
{\it nonzero} $\varphi_k\in\dot P^{\omega_k}_n$ such that
 \begin{equation}\label{Eq phi_1+phi_2+...+phi_n+1=0}
\varphi_1+\dots+\varphi_{n+1}=0
 \end{equation}
holds. Moreover, $\varphi_k$ are unique up to a constant
multiplier: if $\varphi'_k\in\dot P_n^{\omega_k}$ and
$\varphi'_1+\dots+\varphi'_{n+1}=0$ then $\varphi'_k=\lambda
\varphi_k$ with some $\lambda\in\mathbb R$.
\smallskip

\noindent$\bullet$\,\,\,Let $\psi_k\in \dot P_n^{\omega_k}$ be
dual to $\varphi_k$, i.e.,
$\nabla\psi_k=\omega_k\wedge\nabla\varphi_k$. The element
 \begin{equation}\label{Eq psi_1+psi_2+...+psi_n+1=0}
u^*\,=\,\psi_1 \omega_1+\dots+\psi_{n+1}\omega_{n+1}
 \end{equation}
is a {\it nonzero} vector field. Indeed, assuming $u^*=0$, we have
 $$
0=u^*\cdot\omega_{n+1}\,=\,a_1\psi_1+\dots+\psi_{n+1}
 $$
with $a_k=\omega_1\cdot\omega_{n+1}$ and $\psi_k\in \dot
P^{\omega_k}_n$. By the above mentioned uniqueness of the summands
in (\ref{Eq phi_1+phi_2+...+phi_n+1=0}) we get
$\psi_{n+1}=\lambda\varphi_{n+1}$ that contradicts to
$\nabla\psi_{n+1}=\omega_{n+1}\wedge\nabla\varphi_{n+1}$.

As a consequence of (\ref{Eq phi_1+phi_2+...+phi_n+1=0}),\,
(\ref{Eq psi_1+psi_2+...+psi_n+1=0}), we get a {\it nonzero}
quaternion field
 \begin{equation}\label{Eq {0,u*}}
\{0,u^*\}=\{\varphi_1,\psi_1
\omega_1\}+\dots+\{\varphi_{n+1},\psi_{n+1}\omega_{n+1}\}
 \end{equation}
which belongs to $\dot{\mathscr P}_n$, has the zero scalar
component, and is expanded over the axial fields
$\{\varphi_k,\psi_k \omega_k\}\in \dot{\mathscr P}_n^{\omega_k}$.

Since $\{0,u^*\}\in \mathscr P$, one has ${\rm rot\,}u^*=\nabla
0=0$. Hence, $u^*$ is a potential field and one can represent
 \begin{equation}\label{Eq {0,nabla beta*}}
u^*=\nabla \beta^* \quad \text{with} \quad \beta^* \in \dot
P_{n+1}\,.
 \end{equation}

\noindent$\bullet$\,\,\,Let $\dot{\mathscr P}_{0\,n}$ be the
subspace in $\dot{\mathscr P}_n$ of elements with zero scalar
component, so that $\{0,u^*\}\in\dot{\mathscr P}_{0\,n}$. In the
mean time, all elements of $\dot{\mathscr P}_{0\,n}$ are of the
form (\ref{Eq {0,nabla beta*}}), i.e.,
 $$
\dot{\mathscr P}_{0\,n}=\{\{0,\nabla\beta\}\,|\,\,\beta \in \dot
P_{n+1}\}\,.
 $$

By ${\cal O}_3$ we denote the group of rotations of ${\mathbb
R}^3$. As is easy to check, $\dot{\mathscr P}_{0\,n}$ provides a
finite-dimensional representation of the rotation group with the
action
 $$
R:\,\{0,\nabla\beta(x)\}\mapsto\{0,
R\left[\nabla\beta(R^{-1}x)\right]\}\,,\qquad\,\,x\in{\mathbb
R}^3\,\,\,\,(R\in{\cal O}_3)\,.
 $$
As is obvious, such a representation is in fact identical to the
standard representation of ${\cal O}_3$ in $\dot P_{n+1}$. By the
latter, it is irreducible.

In the mean time, the subspace ${\rm span\,}\{
R\{0,u^*\}\,|\,\,R\in{\cal O}_3\}\subset \dot{\mathscr P}_{0\,n}$
is invariant w.r.t. the group action. Hence, the irreducibility
yields
  \begin{equation}\label{Eq *}
{\rm span\,}\{R\{0,u^*\}\,|\,\,R\in{\cal O}_3\}\,=\, \dot{\mathscr
P}_{0\,n}\,.
 \end{equation}

As it easily follows from (\ref{Eq {0,u*}}), the fields
$R\{0,u^*\}$ are also sums of the axial fields. Hence, each
element of the left hand side in (\ref{Eq *}) is a finite sum of
axial fields. Therefore, in accordance with the equality (\ref{Eq
*}), {\it any} element of $\dot{\mathscr P}_{0\,n}$ is a finite
sum of the axial polynomial harmonic fields:
 \begin{equation}\label{Eq span {0,Ru*}}
\dot{\mathscr P}_{0\,n}\,=\,{\rm span}\{\dot{\mathscr
P}^\omega_{0\,n}\,|\,\,\omega \in S^2\}.
 \end{equation}
\smallskip

\noindent$\bullet$\,\,\,Show that
 \begin{equation}\label{Eq Pn=sum Pn omega}
\dot{\mathscr P}_n={\rm span\,}\{\dot{\mathscr
P}_n^\omega\,|\,\,\omega\in S^2\}\,.
 \end{equation}
Take $p=\{\alpha,u\}\in \dot{\mathscr P}_n$. Represent
$\alpha=\phi_1+\dots+\phi_l$ with $\phi_k\in P_n^{\omega_k}$
\footnote{such a representation is not unique but any is
available}. Let $\eta_k\in \dot P_n^{\omega_k}$ satisfy
$\nabla\eta_k=\omega_k\wedge\nabla\phi_k$, so that
$p_k^\prime=\{\phi_k,\eta_k\omega_k\}\in \dot{\mathscr
P}_n^{\omega_k}$. The quaternion polynomial
$p'=p_1^\prime+\dots+p_l^\prime$ belongs to the r.h.s. of (\ref{Eq
Pn=sum Pn omega}). Representing $p=p'+p_0$, one has
$p_0=\{0,u_0\}\in \dot{\mathscr P}_{0\,n}$ by construction. In the
mean time, (\ref{Eq span {0,Ru*}}) yields $p_0$ to be a sum of
axial fields. Hence, eventually, $p$ is also a finite sum of axial
polynomial harmonic fields, i.e., (\ref{Eq Pn=sum Pn omega}) does
hold.
\smallskip

\noindent$\bullet$\,\,\,Notice that the elements of ${\mathscr
P}_0$ \,(constant fields) are axial. Then, representing
 \begin{align*}
& {\mathscr P}={\rm span\,}\{\dot{\mathscr P}_n\,|\,\,n\geqslant
0\}\overset{(\ref{Eq Pn=sum Pn omega})}= {\rm
span\,}\{\dot{\mathscr P}_n^\omega\,|\,\,n\geqslant
0,\,\,\omega\in S^2\}=\\
& =\,{\rm span\,}\{{\mathscr P}^\omega\,|\,\,\omega\in S^2\}\,,
 \end{align*}
we get (\ref{Eq P=sum Pn homog}) and prove Lemma \ref{Lemma P=sum
Pn homog}.

\subsubsection*{Density lemma}
Here we prove Lemma \ref{Lemma on density}.

Take a field $p=\{\alpha,u\}\in{\mathscr Q}(B)$ and show that it
can be approximated  by elements of $\mathscr P$.
\smallskip

\noindent$\bullet$\,\,\,We say $p=\{\alpha,u\}$ to be {\it smooth}
and write $p\in\mathscr S$ if $\alpha\in C^2(B)$ and $u\in C^2(B;
{\mathbb R}^3)$. As is well known, the lineal ${\mathscr
Q}(B)\cap\mathscr S$ is dense in ${\mathscr Q}(B)$.

Let $p=\{\alpha,u\}$ be smooth. Fix a (small) $\varepsilon>0$. For
the harmonic divergence-free field $u$ one can find a harmonic
polynomial vector field $v$ such that $\|u-v\|_{C^2(B;{\mathbb
R}^3)}<\varepsilon$. Since ${\rm div\,}u=0$, the latter inequality
yields
 \begin{equation}\label{Eq estimate div v}
\|{\rm div\,}v\|_{C^1(B)}<{\rm const}\,\varepsilon\,.
 \end{equation}
Also, we have
 \begin{equation}\label{Eq Delta div v=0}
\Delta{\rm div\,}v\,=\,{\rm div\,}\Delta v\,=\,0 \qquad
\text{in}\,\,B\,,
 \end{equation}
so that ${\rm div\,}v$ is a `small' scalar harmonic  polynomial in
the ball.
\smallskip

\noindent$\bullet$\,\,\,By $\partial_{x_k}$ we denote the partial
derivative w.r.t. $x_k$.

Let $q=q(x_1,x_2)$ be a polynomial satisfying
 \begin{equation}\label{Eq Delta q}
\Delta q(x_1, x_2)\,=\,-\partial_{x_3}[{\rm div\,}v](x_1,x_2,0)\,.
 \end{equation}
The function
 $$
\tilde\eta(x_1,x_2,x_3)\,=\,q(x_1,x_2)+\int_0^{x_3}[{\rm
div\,}v](x_1,x_2,t)\,dt
 $$
is a scalar harmonic polynomial. Indeed,
 \begin{align*}
& \Delta \tilde\eta(x_1,x_2,x_3)=\Delta q(x_1,
x_2)+\int_0^{x_3}\left[\partial^2_{x_1}{\rm
div\,}v+\partial^2_{x_2}{\rm div\,}v\right](x_1,x_2,t)\,dt+\\
& +\partial_{x_3}[{\rm div\,}v](x_1,x_2,x_3)\overset{(\ref{Eq Delta div v=0})}=\\
& =\Delta q(x_1, x_2)-\int_0^{x_3}\partial^2_{x_3}[{\rm
div\,}v](x_1,x_2,t)\,dt+\partial_{x_3}[{\rm div\,}v](x_1,x_2,x_3)\overset{(\ref{Eq Delta q})}=\\
& =-\partial_{x_3}[{\rm
div\,}v](x_1,x_2,0)-\left\{\partial_{x_3}[{\rm
div\,}v](x_1,x_2,x_3)-\partial_{x_3}[{\rm
div\,}v](x_1,x_2,0)\right\}+\\
& +\,\partial_{x_3}[{\rm div\,}v](x_1,x_2,x_3)\,=\,0\,.
 \end{align*}

Next, let $r=r(x_1,x_2)$ be a harmonic polynomial satisfying
 \begin{equation*}
\Delta r=0\,, \qquad r=q \quad\text{as}\,\,x_1^2+x_2^2=1\,.
\end{equation*}
By the choice, one has
\begin{equation}\label{Eq Delta q-r}
\Delta (q-r)\,\overset{(\ref{Eq Delta q})}=\,-\partial_{x_3}[{\rm
div\,}v](\cdot,\cdot,0)\,,\qquad q-r=0
\quad\text{as}\,\,x_1^2+x_2^2=1\,.
 \end{equation}
Owing to estimate (\ref{Eq estimate div v}), the well-known
properties of the elliptic Dirichlet problem (\ref{Eq Delta q-r})
provide
 \begin{equation}\label{Eq Estimate q-r}
\|q-r\|_{C^2(B)}\leqslant {\rm const\,}\varepsilon\,.
 \end{equation}

Estimating the integral with regard to (\ref{Eq estimate div v})
and taking into account (\ref{Eq Estimate q-r}), we conclude that
the function $\eta=\tilde\eta-r$,
 $$
\eta(x_1,x_2,x_3)=q(x_1,x_2)-r(x_1,x_2)+\int_0^{x_3}[{\rm
div\,}v](x_1,x_2,t)\,dt
 $$
is a harmonic polynomial, which satisfies
 \begin{equation*}
\partial_{x_3}\eta={\rm div\,}v,\qquad \|\eta\|_{C^1(B)}\leqslant {\rm const\,}\varepsilon\,.
 \end{equation*}

\noindent$\bullet$\,\,\,By the latter, $\eta\bf k$ is a harmonic
polynomial vector field, which satisfies ${\rm div\,}\eta{\bf
k}={\rm div\,}v$ and $\|\eta{\bf k}\|_{C^1(B;{\mathbb R}^3)}<{\rm
const}\,\varepsilon$. As a consequence, $\tilde u=v-\eta{\bf k}$
is a harmonic polynomial vector field satisfying ${\rm
div\,}\tilde u=0$.

Representing $u-\tilde u=u-v-\eta{\bf k}$, we have
 \begin{equation}\label{Eq tilde v approx u}
\|u-\tilde u\|_{C^1(B;{\mathbb
R}^3)}\leqslant\|u-v\|_{C^2(B;{\mathbb R}^3)}+\|\eta{\bf
k}\|_{C^1(B;{\mathbb R}^3)}\leqslant{\rm const}\,\varepsilon\,.
 \end{equation}
So, the vector component $u$ of the smooth harmonic quaternion
field $p$ is approximated by the harmonic polynomial
divergence-free vector field $\tilde u$.

In the mean time, ${\rm rot\,}\tilde u$ is also a harmonic
polynomial divergence-free vector field and the estimate
 \begin{equation}\label{Eq tilde rot tilde u approx rot u}
\|{\rm rot\,}u-{\rm rot\,}\tilde u\|_{C(B;{\mathbb
R}^3)}\leqslant{\rm const}\,\varepsilon
\end{equation}
obviously follows from (\ref{Eq tilde v approx u}).
\smallskip

\noindent$\bullet$\,\,\,By the construction of $\tilde u$, one has
${\rm rot\,rot\,}\tilde u=\nabla{\rm div\,}\tilde u-\Delta\tilde
u=0$. Hence, the field ${\rm rot\,}\tilde u$ has a scalar
potential in $B$. Therefore, integrating over the proper paths in
$B$, one can find a scalar harmonic polynomial $\tilde\alpha$
satisfying $\nabla \tilde\alpha={\rm rot\,}\tilde u$ and
$\tilde\alpha(0,0,0)=\alpha(0,0,0)$.

The potential $\tilde\alpha$ and the vector field $\tilde u$
determine a harmonic polynomial quaternion field $\tilde
p=\{\tilde\alpha,\tilde u\}\in\mathscr P$. The estimate
 $$
\|\nabla \alpha-\nabla \tilde\alpha\|_{C(B;{\mathbb R}^3)}=\|{\rm
rot\,}u-{\rm rot\,}\tilde u\|_{C(B;{\mathbb
R}^3)}\overset{(\ref{Eq tilde rot tilde u approx rot
u})}\leqslant{\rm const}\,\varepsilon
 $$
easily implies $\|\alpha-\tilde\alpha\|_{C(B)}\leqslant{\rm
const}\,\varepsilon$.

Summarizing, we arrive at
 $$
\|p-\tilde p\|_{C^{\mathbb
H}(B)}=\left\{\|\alpha-\tilde\alpha\|^2_{C(B)}+\|u-\tilde
u\|^2_{C(B;{\mathbb R}^3)}\right\}^\frac{1}{2}\leqslant{\rm
const}\,\varepsilon
 $$
and conclude that $\mathscr P$ is dense in the lineal of smooth
fields. Since this lineal is dense in ${\mathscr Q}(B)$, one
obtains
 \begin{equation}\label{Eq overline P=Q(B)}
\overline{\mathscr P}={\mathscr Q}(B)\,.
 \end{equation}
\noindent$\bullet$\,\,\,With regard to ${\mathscr
P}^\omega\subset{\mathscr A}_\omega(B)$, we have
 $$
{\mathscr P}\overset{(\ref{Eq P=sum Pn homog})}={\rm
span\,}\{{\mathscr P}^\omega\,|\,\,\omega\in S^2\}\subset{\rm
span\,}\{{\mathscr A}_\omega(B)\,|\,\,\omega\in S^2\}\,,
 $$
whereas (\ref{Eq overline P=Q(B)}) follows to (\ref{Eq span A
omega = Q(B)}) and proves Lemma \ref{Lemma on density}.

\subsubsection*{Correspondence $\mu\mapsto x^\mu$}
\noindent$\bullet$\,\,\,The quaternion fields of the form
$\pi(x)=\{x\cdot a,Ax\}$, where $a\in{\mathbb R}^3$ is a vector
and $A$ is a constant 3$\times$3-matrix, are called linear. Note
that $\nabla[x\cdot a]=a$,\,\,${\rm rot\,}Ax$ is a constant vector
field, and ${\rm div\,}Ax={\rm tr\,}A$ holds. Therefore, a linear
field is harmonic if and only if $a={\rm rot\,}Ax$ and ${\rm
tr\,}A=0$.

By ${\mathscr L}(B)\subset{\mathscr Q}(B)$ we denote the subspace
of the linear harmonic fields reduced to the ball.
\smallskip

\noindent$\bullet$\,\,\,Let $\omega_1,\, \omega_2,\, \omega_3$ be
a basis in ${\mathbb R}^3$ normalized by
 $$
\omega_k\cdot\omega_l=\delta_{kl};\,\,\,\,\,\omega_1\wedge\omega_2=\omega_3,\,\,\,
\omega_2\wedge\omega_3=\omega_1,\,\,\,\omega_3\wedge\omega_1=\omega_2\,.
 $$
The axial linear fields
 $$
\pi_1(x)=\{x\cdot\omega_1,[x\cdot\omega_2]\omega_3\},\,\pi_2(x)=\{x\cdot\omega_2,[x\cdot\omega_3]\omega_1\},\,
\pi_3(x)=\{x\cdot\omega_3,[x\cdot\omega_1]\omega_2\}
 $$
are called the {\it coordinate fields}. Since
 $$
\nabla [x\cdot\omega_1]=\omega_1=\omega_2 \wedge\omega_3=\nabla
[x\cdot\omega_2] \wedge\omega_3={\rm
rot\,}[x\cdot\omega_2]\omega_3\,,\quad{\rm
div\,}[x\cdot\omega_2]\omega_3=0\,,
 $$
the field $\pi_1$ is harmonic. Analogously, $\pi_2, \pi_3$ are
harmonic.

We omit the proof of the following simple fact: harmonic linear
fields are expanded over the coordinate fields. As example, any
such $\pi$ is uniquely represented in the form
 \begin{equation}\label{Eq pi=a pi 1+b pi 2}
\pi\,=\,{\mathfrak g}\pi_2\,+\,{\mathfrak h}\pi_3
\quad\text{with}\,\,\,{\mathfrak g},{\mathfrak h}\in\mathbb H\,.
 \end{equation}
In particular, the equality
 \begin{equation}\label{Eq pi1= o3 pi2-o2 pi3}
\pi_1\,=\,{\mathfrak o}_3\pi_2\,-\,{\mathfrak o}_2\pi_3
 \end{equation}
holds with ${\mathfrak o}_k=\{0,\omega_k\}\in {\mathbb H}$ and can
be verified by simple calculations.
\smallskip

\noindent$\bullet$\,\,\,Take $\mu \in {\mathscr M}^{\mathbb H}$
and fix $\omega\in S^2$. Since $\mu$ is multiplicative on
${\mathscr A}_\omega(B)$, the image ${\mathbb
A}_\omega=\mu({\mathscr A}_\omega(B))$ is a commutative subalgebra
of $\mathbb H$. By (\ref{Eq l(h)=h}), one has
$\{0,\omega\}\in{\mathbb A}_\omega$. Hence, ${\mathbb
A}_\omega=\{\{a,b\omega\}\,|\,\,a,b\in{\mathbb R}^3\}$ holds.

By the aforesaid, we have \footnote{Here and in what follows, we
regard the linear fields to be reduced on $B$, i.e., regard them
as elements of ${\mathscr L}(B)$.}
 $$
\mu(\pi_1)=\{a_{12},b_{12}\omega_3\},\quad
\mu(\pi_2)=\{a_{23},b_{23}\omega_1\},\quad
\mu(\pi_3)=\{a_{31},b_{31}\omega_2\}
 $$
with some $a_{kl},b_{kl}\in\mathbb R$.

Applying $\mu$ to (\ref{Eq pi1= o3 pi2-o2 pi3}), we get
 $$
\mu(\pi_1)={\mathfrak o}_3 \mu(\pi_2)\,-\,{\mathfrak o}_2
\mu(\pi_3)
 $$
or, in detail,
 \begin{align*}
&
\{a_{12},b_{12}\omega_3\}=\{0,\omega_3\}\{a_{23},b_{23}\omega_1\}-\{0,\omega_2\}\{a_{31},b_{31}\omega_2\}=\\
& = \{0,a_{23}\omega_3+b_{23}\omega_2\}-\{-b_{31},
a_{31}\omega_2\}=\{b_{31},a_{23}\omega_3+[b_{23}-a_{31}]\omega_2\}\,.
 \end{align*}
Comparing by components, we get
$a_{12}=b_{31},\,\,b_{12}=a_{23},\,\,b_{23}=a_{31}$.
\smallskip

\noindent$\bullet$\,\,\,Denote
$a=a_{12}=b_{31},\,\,b=b_{12}=a_{23},\,\,c=b_{23}=a_{31}$ and
define a vector
 $$
x^\mu\,=\,a\omega_1+b\omega_2+c\omega_3\,\in{\mathbb R}^3\,.
 $$
Summarizing the previous calculations, we easily conclude that
 \begin{equation}\label{Eq point x^mu}
\mu(\pi_1)=\pi_1(x^\mu),\,\,\,\,\mu(\pi_2)=\pi_2(x^\mu),\,\,\,\,\mu(\pi_3)=\pi_3(x^\mu)\,.
 \end{equation}
So, the functional $\mu$,which is $\mathbb H$-linear and
multiplicative on the axial algebras, determines a point
$x^\mu\in{\mathbb R}^3$ associated with it via (\ref{Eq point
x^mu}).

\subsubsection*{Completing proof of Theorem \ref{T2}}

\noindent$\bullet$\,\,\,Applying $\mu$ to (\ref{Eq pi=a pi 1+b pi
2}), we have
 \begin{align}
\notag & \mu(\pi)=\mu({\mathfrak g}\pi_2\,+\,{\mathfrak
h}\pi_3)={\mathfrak g}\mu(\pi_2)\,+\,{\mathfrak
h}\mu(\pi_3)\overset{(\ref{Eq point x^mu})}={\mathfrak
g}\pi_2(x^\mu)\,+\,{\mathfrak
h}\pi_3(x^\mu)\,=\\
& = \,\left[{\mathfrak g}\pi_2\,+\,{\mathfrak
h}\pi_3\right](x^\mu)\,=\,\pi(x^\mu)\,\label{Eq final}
 \end{align}
that extends (\ref{Eq point x^mu}) to ${\mathscr L}(B)$.
\smallskip

\noindent$\bullet$\,\,\,As is well known, the disk algebra
${\mathscr A}(D)$ is generated by the functions $1,\,z$. By
perfect analogy with this fact, any axial algebra ${\mathscr
A}_\omega(B)\cong{\mathscr A}(D)$ is generated by two fields
$\{1,0\}$ and $\{x\cdot \eta,[x\cdot
\eta\wedge\omega]\omega\}\in{\mathscr L}(B)$, where $\eta\in S^2$
is arbitrary provided $\eta\cdot\omega=0$. Therefore, (\ref{Eq
final}) is extended to ${\mathscr A}_\omega(B)$ by continuity and
yields
 $$
\mu(p)\,=\,p(x^\mu)\,,\qquad p\in {\mathscr
A}_\omega(B)\qquad(\omega\in S^2)\,.
 $$
Here we see that the point $x^\mu$ must belong to the ball $B$.
Otherwise, one can chose (for instance, the polynomials) $p_j\in
{\mathscr A}_\omega(B)$ such that $\|p_j\|_{{\mathscr
Q}(B)}\leqslant \rm const$ and $|p_j(x^\mu)|\to \infty$ in
contradiction to the continuity of $\mu$.

Then, in accordance with Lemma \ref{Lemma on density}, one extends
$\mu$ to ${\mathscr Q}(B)$ an obtains $\mu(p)\,=\,p(x^\mu)$ for
all $p\in {\mathscr Q}(B)$. Hence, we conclude that
$\mu=\delta^{\mathbb H}_{x^\mu}$ is valid.
\smallskip

\noindent$\bullet$\,\,\,The latter equality easily implies that
the bijection ${\mathscr M}^{\mathbb H}\ni\mu\mapsto x^\mu\in B$
is a homeomorphism of topological spaces. To show this one can
just check that $\mu_j\to\mu$ in ${\mathscr M}^{\mathbb H}$ is
equivalent to $x^{\mu_j}\to x^\mu$ in $B$.
\smallskip

The proof of Theorem \ref{T2} is completed.

\end{document}